\documentclass[12pt,twoside,a4paper]{article}
\usepackage[english]{babel}
\usepackage{amssymb}
\usepackage[T1]{fontenc}
\usepackage{amsmath,epsf}
\newtheorem{Th}{Theorem}
\newtheorem{lemma}{Lemma}
\newtheorem{Cor}{Corollary}

\newcommand{\A}{\mathcal{A}}

\renewcommand{\L}{\mathcal{L}}
\newcommand{\F}{\mathcal{F}}
\newcommand{\R}{\mathbb{R}}
\newcommand{\Z}{\mathbb{Z}}

\renewcommand{\P}{\mathbb{P}}

\newcounter{tictac}

\def\1{\,\rlap{\mbox{\small\rm 1}}\kern.15em 1}

\def\build#1_#2^#3{\mathrel{\mathop{\kern 0pt#1}\limits_{#2}^{#3}}}
\def\tend#1#2{\build\hbox to 12mm{\rightarrowfill}_{#1\rightarrow #2}^{a.s.}}

\def\converge#1#2#3{\build\hbox to
15mm{\rightarrowfill}_{#1\rightarrow #2}^{\hbox{\scriptsize #3}}}
\begin{document}
\title{Berry-Esseen's central limit theorem for non-causal linear processes in Hilbert space}
\author{Mohamed EL MACHKOURI}
\maketitle

{\renewcommand\abstractname{Abstract}
\begin{abstract}
\baselineskip=18pt 
Let $H$ be a real separable Hilbert space and $(a_k)_{k\in\Z}$ a sequence of bounded linear operators from $H$ to $H$. 
We consider the linear process $X$ defined for any $k$ in $\Z$ by $X_k=\sum_{j\in\Z}a_j(\varepsilon_{k-j})$ where 
$(\varepsilon_k)_{k\in\Z}$ is a sequence of i.i.d. centered $H$-valued random variables. We investigate the rate of convergence 
in the CLT for $X$ and in particular we obtain the usual Berry-Esseen's bound provided that 
$\sum_{j\in\Z}\vert j\vert\|a_j\|_{\L(H)}<+\infty$ and $\varepsilon_0$ belongs to $L_H^{\infty}$.\\
\\
{\em Short title:} Berry-Esseen's CLT for Hilbertian linear processes.\\
{\em Key words:} Central limit theorem, Berry-Esseen bound, linear process, Hilbert space.
\end{abstract}
\thispagestyle{empty}
\baselineskip=18pt
\section{Introduction and notations}
Let $(H,\|.\|_H)$ be a separable real Hilbert space and $(\L,\|.\|_{\L(H)})$ be the class of bounded linear 
operators from $H$ to $H$ with its usual uniform norm. Consider a sequence $(\varepsilon_k)_{k\in\Z}$ of i.i.d. centered 
random variables, defined on a probability space $(\Omega,\A,\P)$, with values in $H$. 
If $(a_k)_{k\in\Z}$ is a sequence in $\L$, we define the (non-causal) linear process $X=(X_k)_{k\in\Z}$ in $H$ by 
\begin{equation}\label{definition_linear_process}
X_k=\sum_{j\in\Z}a_j(\varepsilon_{k-j}),\qquad k\in\Z.
\end{equation}
If $\sum_{j\in\Z}\|a_j\|_{\L(H)}<\infty$ and $E\|\varepsilon_0\|_H<+\infty$ then the series in 
$(\ref{definition_linear_process})$ converges almost surely and in $L_H^1(\Omega,\A,\P)$ (see Bosq \cite{Bosq___2000}). 
The condition $\sum_{j\in\Z}\|a_j\|_{\L(H)}<\infty$ is know to be sharp for the $\sqrt{n}$-normalized 
partial sums of $X$ to satisfies a CLT provided that $(\varepsilon_k)_{k\in\Z}$ are i.i.d. centered having finite 
second moments (see Merlevede et al. \cite{Merlevede--Peligrad--Utev___1997}). In this work, we investigate the rate of convergence 
in the CLT for $X$ under the condition
\begin{equation}\label{condition}
\sum_{j\in\Z}\vert j\vert^{\tau}\|a_j\|_{\L(H)}<\infty
\end{equation}
with $\tau=1$ when $(\varepsilon_k)_{k\in\Z}$ are assumed to be i.i.d. centered and such that $\varepsilon_0$ belongs to 
$L_H^{\infty}$ and $\tau=1/2$ when $(\varepsilon_k)_{k\in\Z}$ are i.i.d. centered and such that $\varepsilon_0$ belongs 
to some Orlicz space $L_{H,\psi}$ (see section 2). 
This problem was previously studied (with $\tau=1$ in Condition $(\ref{condition})$) by 
Bosq \cite{Bosq___2003} for (causal) Hilbert linear processes but a 
mistake in his proof was pointed out by V. Paulauskas \cite{Paulauskas___2004}. However, in the particular case of Hilbertian autoregressive processes 
of order 1, Bosq \cite{Erratum_Bosq___2004} obtained the usual Berry-Esseen inequality provided that 
$(\varepsilon_k)_{k\in\Z}$ are i.i.d. centered with $\varepsilon_0$ in $L_H^{\infty}$.
\section{Main result}
In the sequel, $C_{\varepsilon_0}$ is the autocovariance operator of $\varepsilon_0$, $A:=\sum_{j\in\Z}a_j$ and $A^{\ast}$ 
is the adjoint of $A$. For any sequence $Z=(Z_k)_{k\in\Z}$ of random variables with values in $H$ we denote
$$
\Delta_n(Z)=\sup_{t\in\R}\bigg\vert\P\left(\bigg\|\frac{1}{\sqrt{n}}\sum_{k=1}^nZ_k\bigg\|_H\leq t\right)
-\P\left(\|N\|_H\leq t\right)\bigg\vert
$$
where $N\sim\mathcal{N}(0,AC_{\varepsilon_0}A^{\ast})$.\\
\\
For any $j\in\Z$, denote $c_{j,n}=\sum_{i=1}^nb_{i-j}$ where $b_i=a_i$ for any $i\neq 0$ and $b_0=a_0-A$.
\begin{lemma}\label{approximation-sommes-partielles}
For any positive integer $n$,
$$
\sum_{k=1}^nX_k=A\left(\sum_{k=1}^n\varepsilon_k\right)+Q_n+R_n
$$
where $Q_n=\sum_{k=1}^n\sum_{\vert j\vert>n}a_{k-j}(\varepsilon_{j})$ and $R_n=\sum_{\vert j\vert\leq n}c_{j,n}(\varepsilon_j)$.
\end{lemma}
Recall that a Young function $\psi$ is a real convex nondecreasing
function defined on $\R^{+}$ which satisfies $\lim_{t\to+\infty}\psi(t)=+\infty$ and $\psi(0)=0$.
We define the Orlicz space $L_{H,\psi}$ as the space of $H$-valued random
variables $Z$ defined on the probability space $(\Omega, \F, \P)$
such that $E[\psi(\| Z\|_H/c)]<+\infty$ for some $c>0$. The
Orlicz space $L_{H,\psi}$ equipped with the so-called Luxemburg
norm $\| . \|_{\psi}$ defined for any $H$-valued random
variable $Z$ by
$$
\| Z\|_{\psi}=\inf\{\,c>0\,;\,E[\psi(\|Z\|_{H}/c)]\leq 1\,\}
$$
is a Banach space. In the sequel, $c(N)$ denotes a bound of the density of 
$\mathcal{N}(0,AC_{\varepsilon_0}A^{\ast})$ (see Davydov et al. \cite{Davydov--Lifshits--Smorodina___1998}). 
Our main result is the following. 
\begin{Th}\label{mainresult} Let $(\varepsilon_k)_{k\in\Z}$ be a sequence of i.i.d. centered $H$-valued random variables 
and let $X$ be the Hilbertian linear process defined by $(\ref{definition_linear_process})$. 
\begin{itemize}
 \item[i)] If $\varepsilon_0$ belongs to $L_H^{\infty}$ and 
$\sum_{j\in\Z}\vert j\vert \| a_j\|_{\L(H)}<\infty$ then 
\begin{equation}\label{BE-inequality_L-infinity}
\Delta_n(X)\leq\frac{c_1}{\sqrt{n}} 
\end{equation}
where $c_1=c_2+14c(N)\|\varepsilon_0\|_{\infty}\sum_{j\in\Z}\vert j\vert \| a_j\|_{\L(H)}$ and $c_2$ is a positive 
constant which depend only on the distribution of $\varepsilon_0$.
 \item[ii)] If $\psi$ is a Young function then
\begin{equation}\label{BE-inequality_L-psi}
\Delta_n(X)\leq\Delta_n(A(\varepsilon))+\varphi\left(\frac{c(N)\|Q_n+R_n\|_{\psi}}{\sqrt{n}}\right)
\end{equation}
where $\varphi(x)=xh^{-1}(1/x)$ and $h(x)=x\psi(x)$ for any real $x>0$. 
\end{itemize}
\end{Th}
The inequality ($\ref{BE-inequality_L-psi}$) ensures a rate of convergence to zero for 
$\Delta_n(X)$ as $n$ goes to infinity provided that $\Delta_n(A(\varepsilon_0))$ goes to zero as $n$ goes to infinity and 
a bound for $\|Q_n+R_n\|_{\psi}$ exists. As just an illustration, we have the following corollary.
\begin{Cor}\label{corollary-mainresult}
Assume that $(\varepsilon_k)_{k\in\Z}$ are i.i.d. centered $H$-valued random variables and that the condition 
$(\ref{condition})$ holds with $\tau=1/2$.
\begin{itemize}
\item[i)] If $\varepsilon_0$ belongs to $L_{H,\psi_1}$ then 
$\Delta_n(X)=O\left(\frac{\log n}{\sqrt{n}}\right)$ where $\psi_1$ is the Young function defined 
by $\psi_1(x)=\exp(x)-1$.
\item[ii)] If $\varepsilon_0$ belongs to $L_H^r$ for $r\geq 3$ then 
$\Delta_n(X)=O\left(n^{-\frac{r}{2(r+1)}}\right)$.
\end{itemize}
\end{Cor}

\section{Proofs}
{\bf Proof of Lemma $\textbf{\ref{approximation-sommes-partielles}}$}. For any positive integer $n$, we have
\begin{align*}
R_n&=\sum_{j=-n}^nc_{j,n}(\varepsilon_j)=\sum_{k=1}^n\sum_{j=-n}^nb_{k-j}(\varepsilon_j)\\
&=\sum_{k=1}^n\sum_{j\in [-n,n]\backslash\{k\}}a_{k-j}(\varepsilon_j)+(a_0-A)\left(\sum_{k=1}^n\varepsilon_k\right)\\
&=\sum_{k=1}^n\sum_{j=-n}^na_{k-j}(\varepsilon_j)-A\left(\sum_{k=1}^n\varepsilon_k\right)\\
&=-Q_n+\sum_{k=1}^n X_k-A\left(\sum_{k=1}^n\varepsilon_k\right).
\end{align*}
The proof of Lemma $\ref{approximation-sommes-partielles}$ is complete.\\
\\
%
%
{\bf Proof of Theorem $\textbf{\ref{mainresult}}$}. Let $\lambda>0$ and $t>0$ be fixed and denote 
$U=A\left(\sum_{k=1}^n\varepsilon_k/\sqrt{n}\right)$ and $V=(Q_n+R_n)/\sqrt{n}$. So $U+V=\sum_{k=1}^nX_k/\sqrt{n}$ and 
\begin{equation}\label{crucial-inequality}
\P(\|U+V\|_H\leq t)\leq\P(\|U\|_H\leq t+\lambda)+\P(\|V\|_H\geq\lambda)
\end{equation}
For $\lambda_0=2\|V\|_{\infty}$, we obtain
$$
\P(\|U+V\|_H\leq t)-\P(\|N\|_H\leq t)
\leq\P(\|U\|_H\leq t+\lambda_0)-\P(\|N\|_H\leq t).
$$
If $c(N)$ denotes a bound for the density of $\|N\|_H$ (see Davydov et al. \cite{Davydov--Lifshits--Smorodina___1998}) then 
$$
\Delta_n(X)\leq\Delta_n(A(\varepsilon))+\frac{2c(N)\|Q_n+R_n\|_{\infty}}{\sqrt{n}}.
$$
Noting that 
\begin{equation}\label{decomposition-Qn}
Q_n=\sum_{j\geq n+2}a_j\left(\sum_{k=1-j}^{-n-1}\varepsilon_k\right)+\sum_{j<0}a_j\left(\sum_{k=n+1}^{n-j}\varepsilon_k\right)
\end{equation}
and 
\begin{equation}\label{decomposition-Rn}
R_n=R_n^{'}+R_n^{''} 
\end{equation}
where
$$
R_n^{'}=-\sum_{j=-n}^{-1}a_j\left(\sum_{k=1}^{-j}\varepsilon_k\right)-\sum_{j<-n}a_j\left(\sum_{k=1}^n\varepsilon_k\right)-\sum_{j>0}a_j\left(\sum_{k=n-j+1}^n\varepsilon_k\right)\\
$$
and
$$
R_n^{''}=\sum_{j=1}^na_j\left(\sum_{k=-j+1}^0\varepsilon_k\right)+\sum_{j=n+1}^{2n}a_j\left(\sum_{k=-n}^{n-j}\varepsilon_k\right),
$$
we derive that $\|Q_n+R_n\|_{\infty}\leq 7\|\varepsilon_0\|_{\infty}\sum_{j\in\Z}\vert j\vert\|a_j\|_{\L(H)}$ and consequently
$$
\Delta_n(X)\leq\Delta_n(A(\varepsilon))+\frac{14c(N)\|\varepsilon_0\|_{\infty}\sum_{j\in\Z}\vert j\vert\|a_j\|_{\L(H)}}{\sqrt{n}}.
$$
Combining the last inequality with the Berry-Esseen inequality for i.i.d. centered $H$-valued random variables 
(see Yurinski \cite{Yurinskii___1982} or Bosq \cite{Bosq___2000}, Theorem 2.9) we obtain ($\ref{BE-inequality_L-infinity}$).\\
\\
In the other part, if $\psi$ is a Young function we have 
$\P(\|V\|_H\geq\lambda)\leq\frac{1}{\psi\left(\lambda/\|V\|_{\psi}\right)}$ 
and keeping in mind inequality $(\ref{crucial-inequality})$, we derive
$$
\Delta_n(X)\leq\Delta_n(A(\varepsilon))+c(N)\lambda+\frac{1}{\psi\left(\lambda/\|V\|_{\psi}\right)}.
$$
Noting that 
$c(N)\lambda=\frac{1}{\psi\left(\lambda/\|V\|_{\psi}\right)}$ if and only if 
$\lambda=\frac{\varphi\left(c(N)\|V\|_{\psi}\right)}{c(N)}$ where $\varphi$ is defined by
$\varphi(x)=xh^{-1}(1/x)$ and $h$ by $h(x)=x\psi(x)$, we conclude
$$
\Delta_n(X)\leq\Delta_n(A(\varepsilon))+\varphi\left(\frac{c(N)\|Q_n+R_n\|_{\psi}}{\sqrt{n}}\right).
$$
The proof of Theorem $\ref{mainresult}$ is complete.\\
\\
{\bf Proof of Corollary $\textbf{\ref{corollary-mainresult}}$.} Assume that $\|\varepsilon_0\|_{\psi_1}<\infty$ 
where $\psi_1$ is the Young function defined by $\psi_1(x)=\exp(x)-1$. There exists $a>0$ such that 
$E\left(\exp(a\|\varepsilon_0\|_H)\right)\leq 2$. So, there exist 
(see Arak and Zaizsev \cite{Arak--Zaitsev___1988}) constants $B$ and $L$ such that 
$$
E\|\varepsilon_0\|_H^{m}\leq\frac{m!}{2}B^2L^{m-2},\quad m=2,3,4,...
$$
Applying Pinelis-Sakhanenko inequality 
(see Pinelis and Sakhanenko \cite{Pinelis-Sakhanenko___1986} or Bosq \cite{Bosq___2000}), we obtain
$$
\P\left(\bigg\|\sum_{k=p}^q\varepsilon_k\bigg\|_H\geq x\right)
\leq\exp\left(-\frac{x^2}{2(q-p+1)B^2+2xL}\right),\quad x>0
$$
and using Lemma 2.2.10 in Van Der Vaart and Wellner \cite{an-Der-Vaart--Wellner___1996}, there exists a universal 
constant $K$ such that
\begin{equation}\label{psi1-norm}
\bigg\|\sum_{k=p}^q\varepsilon_k\bigg\|_{\psi_1}\leq K\left(L+B\sqrt{q-p+1}\right)
\end{equation}
Combining ($\ref{decomposition-Qn}$), ($\ref{decomposition-Rn}$) and ($\ref{psi1-norm}$), we derive 
$\|Q_n+R_n\|_{\psi_1}\leq C\sum_{j\in\Z}\sqrt{\vert j\vert}\|a_j\|_{\L(H)}$ where the constant $C$ does not depend on $n$. 
Keeping in mind the Berry-Esseen's central limit theorem for i.i.d. centered $H$-valued random variables 
(see Yurinski \cite{Yurinskii___1982} or Bosq \cite{Bosq___2000}, Theorem 2.9), we apply Theorem $\ref{mainresult}$ with the 
Young function $\psi_1$. Since the function $\varphi$ defined by $\varphi(x)=xh^{-1}(1/x)$ with $h(x)=x\psi_1(x)$ satisfies 
$$
\lim_{x\to 0}\frac{\varphi(x)}{x\log(1+\frac{1}{x})}=0,
$$
we derive $\Delta_n(X)=O\left(\frac{\log n}{\sqrt{n}}\right)$.\\
\\
Now, assume that $\|\varepsilon_0\|_{r}<\infty$ for some $r\geq 3$. Applying Pinelis inequality 
(see Pinelis \cite{Pinelis___1994}), there exists a universal constant $K$ such that
$$
\bigg\|\sum_{k=p}^q\varepsilon_k\bigg\|_r\leq K\left(r\left(\sum_{k=p}^q E\|\varepsilon_k\|_H^r\right)^{1/r}+\sqrt{r}\left(\sum_{k=p}^q E\|\varepsilon_k\|_H^2\right)^{1/2}\right)
$$
and consequently
\begin{equation}\label{Lp-norm}
\bigg\|\sum_{k=p}^q\varepsilon_k\bigg\|_r\leq 2Kr\|\varepsilon_0\|_{r}\sqrt{q-p+1}.
\end{equation}
Combining ($\ref{decomposition-Qn}$), ($\ref{decomposition-Rn}$) and ($\ref{Lp-norm}$), we derive $\|Q_n+R_n\|_{r}\leq C\sum_{j\in\Z}\sqrt{\vert j\vert}\|a_j\|_{\L(H)}$ 
where the constant $C$ does not depend on $n$. Again, applying Berry-Esseen's CLT (see Yurinski \cite{Yurinskii___1982} or Bosq \cite{Bosq___2000}, Theorem 2.9) and Theorem $\ref{mainresult}$ with the Young 
function $\psi(x)=x^r$ and the function $\varphi$ given by $\varphi(x)=x^{r/(r+1)}$, we obtain 
$\Delta_n(X)=O\left(n^{-\frac{r}{2(r+1)}}\right)$. The proof of Corollary $\ref{corollary-mainresult}$ is complete.


\vspace{+0.1cm}
\noindent
Mohamed EL MACHKOURI\\
Laboratoire de Math\'ematiques Rapha\"el Salem\\
UMR CNRS 6085, Universit\'e de Rouen\\
Avenue de l'universit\'e\\
76801 Saint-Etienne du Rouvray\\
mohamed.elmachkouri@univ-rouen.fr
\end{document}